\documentclass[11pt]{amsart}

\usepackage{amsmath,amsthm,amscd,euscript}

\newtheorem{theorem}{Theorem}[section]
\newtheorem{lemma}[theorem]{Lemma}

\theoremstyle{remark}

\theoremstyle{definition}

\title{Bounds for solution of linear diophantine  equations}
\author{S.~I.~Veselov}
\date{}

\thanks{E-mail: vesi@uic.nnov.ru.}
\begin{document}
\maketitle

{\bf Abstract.}
{\small Given $A\in \mathbf{Z}^{m\times n},\ rank A=m,\ b\in \mathbf{Z}^{m}$. Let $d$ be the maximum of absolute values of the $m\times m$ minors of the matrix $(A\: b)$, $M=\{x\in \mathbf{Z}^n | Ax=b,\ x\ge 0\}.$ It is shown that if $M\not=\emptyset$, then there exists
$x^0=(x^0_1,...,x^0_n)\in M,$ such that $x^0_i\le d\ (i=1,2,..,n).$
}
\bigskip

\centerline{\bf Introduction}
\medskip

Let $A(m\times n)$ be a matrix of rank m with integer elements, $b$ be an integer
vector, $d$ be the maximum of the absolute values of the $m\times m$ minors of $(A\: b)$, $M$ be the set of nonnegative integer solutions for the system
$Ax=b$, $N=\{1,...,n\}$.

In $[1]$ the conjecture was made that the following theorem is true:

\begin{theorem} If $M$ is nonempty, then there exists  $x^0\in M$ such that
$x_i^0\le d, i\in N.$ \end{theorem}
This conjecture was considered in [1-3], howerever, full and strict
answer had not given. The complete proof of the theorem was given in [4]. See also [5]. 
In this paper I state the translation of the proof from [4].

\paragraph{Notation.} Let $H$ denotes the matrix
which rows are the lattice basis for integer solutions of system
$Ax=0$; 

$h_{1},...,h_{n}$ be columns of $H$;

 $a_{1},...,a_{n}$ be columns of $A$; 
 
 $x^{1}$ be any vector of $M$;
 
 $M^{1}$ be the set of integer solutions of the system
 $H^{T}y+x^{1}\ge0$.
 \medskip
 
  \centerline{\bf The proof of Theorem 1}
 \medskip
 
 Without loss of generality, assume that g.c.d of all $m\times m$
 minors of $A$ is equal to $1$. It is clear that the relation
 $x=H^{T}y+x^{1}$  determines one-to-one mapping between $M$ and $M^{1}$.
 We should use next result from [6]:

\begin{lemma}. Let $I\subseteq N,\mid I \mid=m,A^{1}$ be the matrix that consists
 the columns of matrix $H$ with index from $N\setminus I$. Then
 $\mid\det A^{1}\mid=\mid\det H^{1}\mid$ .\end{lemma}

The theorem can be proved by induction on $n$. Case $n=m$ 
 is obvious. 
  Assume that theorem is true for
${n\le n_0}$ and prove it is true for $n=n_0+1.$

1. Suppose first that from ${u^TA\ge0}$ it follows that
$u=0.$  According to Minkowsky-Farkas' theorem the cone
$$Ax=0,x\ge 0$$ has dimension $n-m.$ Then cone $H^Ty\ge 0$
has dimension $n-m$ as well. Suppose without loss of generality
that equation $h_1^Ty=0$ determines $(n-m-1)$  dimensional
face of cone and $x_1^1 =max(-h_1^Ty)$ for subject
$H^Ty+x^1\ge 0.$ Let $s$ be g.c.d. of components of vector
$h_1$ and $r$ be minimal nonnegative integer number such
that $x_1^1-r$ is divided to $s$. Then the set $M^1$ is
described by system
\begin{center}
   \parbox{10cm}
   {
   $h_1^Ty+x_1^1-r\ge 0\\\rm h_i^Ty+x_i^1\ge 0,
       i\in \{2,3,..,n\}
      $
   }
\end{center}
and there is $y^0\in M^1$, such that $h_1^Ty^0+x_1^1-r= 0$, hence the
system $a_2x_2+...+a_nx_n=b-ra_1$ has a nonnegative integer solution.
From the lemma it follows that $s$ is divisor of each minor of extended
matrix of the system. Therefore there exists a matrix $D$
with determinant $s$
such that vectors $ a_i^1=D^{-1}a_i, i\in\{ 2,3,..,n\}$ and
$b^1=D^{-1}(b-ra_1)$ have integer components. Since $r\le s-1$ ,maximal
absolute value of minors with rank $m$ is not more than d and by
induction the system has solution
$(x_2^0,x_3^0,...,x_n^0)$,
which components  do not exceed    $ d$. As $x^0$ it is
possible to choose $(r,x_2^0,x_3^0,...,x_n^0)$

2. Suppose now that there exists $u\ne 0$ such that $u^TA\ge 0.$
If $u^TA > 0$ then $M$ is bounded, hence,  the theorem is true.
Otherwise, we can do an unimodular transformation for rows of matrix
$(A b)$, permute columns and get the matrix
$\left(\begin{array}{l}\rm
a_1^1...a_v^1\;a_{v+1}^1...a_n^1\;b^1\\\medskip\rm
0\;...\;0\;\;a_{v+1}^2...a_n^2\;b^2\end{array}\right)$
where submatrix $(a_1^1...a_v^1)$ has rank $k$ and consists of
$k$ rows
and columns $(a_{v+1}^2...a_n^2)$ have positive last component.
Set $ b^3=b^1-a^1_{v+1}x^1_{v+1}-...--a^1_{n}x^1_{n}$ and consider
the system \begin{equation} a^1_1x_1^1+....+a^1_vx_v^1=b^3 .\label{eq1}\end{equation}
We will prove that if $d^1$ is maximal absolute value of $k\times k$ minors of this
system, then $d^1$ does not exceed $d.$

Let $k\ge 2$. Suppose without loss of generality that $d^1=abs(det(
a_1^1\;...a_{k-1}\; b^3)).$ Show that it is possible to choose the set
$I=\{i_1,...,i_j\}$ , where $j=n-m$, so that determinants

$$det\left(\begin{array}{l}\rm a_1^1...a_{k-1}^1\;a_{i_1}^1...a_{i_j}^1\;b^3
\smallskip
\\
\rm 0\;...\;0\;\;\;\;\;a_{i_1}^2...a_{i_j}^2\;0\end{array}\right)
,
det\left(\begin{array}{l}
\rm a_1^1...a_{k-1}^1\;a_{i_1}^1...a_{i_j}^1\;a_i^1\smallskip      \\\rm
0\;...\;0\;\;\;\;\;a_{i_1}^2...a_{i_j}^2\;a_i^2\end{array}\right),i\notin I$$

are either nonpositive or nonnegative.

Define $\lambda\ne0$ so that $\lambda^ta^1_i,i\in \{1,...,k-1\}$ and
$\lambda^Tb^3>0.$
Further find a vertex $\mu=(\mu_1,...,\mu_n)$ of the set of solution  of the
system $ u^Ta_i^2+\lambda^Ta_i^1\ge 0, i\in\{v+1,...,n\}.$
As the set $I$ we
choose the set of numbers of linear indepedent inequalities wich are equalities
on $\mu.$  Let $\lambda_e\ne 0$ ,then replace $e$-th row of each
determinant by linear combination of rows with coefficients
$\lambda_1,...,\lambda_k,\mu_1,...,\mu_{m-k}$.
Decomposing each determinant on the $e$-th row we get the sequence of numbers
$ Z\lambda_e\lambda^T b^3, Z\lambda_e(\lambda^Ta_i^1+\mu^Ta_i^2),i\notin I$
where $Z$ is value of algebraic supplement. Obviosly, all numbers
are either nonpositive or nonnegative. 

Now we have
$$d\ge
| det\left(\begin{array}{l}\rm a_1^1...a_{k-1}^1\;a_{i_1}^1...a_{i_j}^1\;b^1
\smallskip      \\\rm
0\;...\;0\;\;\;\;\;a_{i_1}^2...a_{i_j}^2\;b^2\end{array}\right)|=
|\sum_{i\in
I}det\left(\begin{array}{l}\rm a_1^1...a_{k-1}^1\;a_{i_1}^1...a_{i_j}^1\;a_i^1
\smallskip      \\\rm
0\;...\;0\;\;\;\;\;a_{i_1}^2...a_{i_j}^2\;a_i^2\end{array}\right)x_i^1|=$$
$$|det\left(\begin{array}{l}
\rm a_1^1...a_{k-1}^1\;a_{i_1}^1...a_{i_j}^1\;b^3\smallskip      \\\rm
0\;...\;0\;\;\;\;\;a_{i_1}^2...a_{i_j}^2\;0\end{array}\right)x_i^1|\ge
|det\left(\begin{array}{l}\rm a_1^1...a_{k-1}^1\;a_{i_1}^1...a_{i_j}^1\;b^3
\smallskip      \\\rm
0\;...\;0\;\;\;\;\;a_{i_1}^2...a_{i_j}^2\;0\end{array}\right)x_i^1|\ge d^1$$

If $k=1$, then put
$\lambda_1=\left\{ \begin{array}{l}\rm\; 1\; for\; b^3>0,
\smallskip      \\\rm -1\; for\; b^3\le 0,\;e=1\end{array}\right.$
and repeat reasoning.

By induction hypothesis the system (\ref{eq1}) has solution
$(x^0_1,...,x^0_v)$, and $x^0_i\le d$ for $ i\in \{1,2,...,v\}.$
Since $x_i^1$ for $i\in\{v+1,...,n\}$ is bounded from above by $\max\{ x_i|\,Ax=b,\;x\ge 0\},$ hence,
it does not exceed $d$. So we may take $x^0=(x^0_1,...,x^0_v,x^1_{v+1},...,x^1_n).$ The proof is completed.

A little modification of the above proof allow to prove that there exists vertex $x$ of $conv(M)$ such that 
 $x_i\le d$ for all $i\in N.$
\begin{center} ╤яшёюъ ышЄхЁрЄєЁ√.
\end{center}

1. Borosh I.,Treybig B. Bounds on positive integral solutions of linear
diophantine equations.// Proc. Amer. Math. Soc. - 1976, v.55,N 2, pp.299-304.

2. Borosh I., Flahive M., Rubin D. and Treybig B. A sharp Bound
for  solutions of linear diophantine equations.// Proc. Amer. Math. Soc. -
1989, v.105,N 4, pp.844-846.

3. Ivanov N.N. Upper bounds for solution of integer programming. // Kibernetika, (6):112-114,1988.

4. Veselov S. I. The proof of one conjecture linear diophantine equations. (Russian) Manuscript
667ЦB93, deposited at VINITI, Moscow, 1993.

5. Veselov, S.I. The proof of a generalization of Borosh-Treybig's hypothesis for Diophantine equations. (Russian)
Diskretn. Anal. Issled. Oper., Ser. 1 8, No.1, 17-22 (2001).

6. V. N. Shevchenko. Qualitative Topics in Integer Linear Programming,
volume 156 of Translations of Mathematical Monographs.
AMS, Providens, Rhode Island, 1997.

\end{document}